\documentclass[12pt,a4paper,reqno,twoside]{amsart}
\usepackage{amsfonts, amsthm, amsmath, amssymb}

\usepackage{hyperref}
\hypersetup{colorlinks=false}

\usepackage[margin=1in]{geometry}
\usepackage[mathscr]{eucal}
\usepackage{helvet}
  
  \usepackage{setspace}
  \usepackage{appendix}
  
\usepackage{amsthm}
\newtheorem{theorem}{Theorem}
\newtheorem{lemma}{Lemma}[section]
\newtheorem{remark}{Remark}

\newtheorem{proposition}{Proposition}
\newtheorem{corollary}{Corollary}
\newcommand{\di}{\digamma}
\newcommand{\ve}{\varepsilon}
\newcommand{\bm}{\boldsymbol}
 \newcommand\blfootnote[1]{
	\begingroup
	\renewcommand\thefootnote{}\footnote{#1} 
	\endgroup
}

\newcommand{\ot}{\otimes}

\newcommand{\bmu}{\bm\mu}
\newcommand{\bnu}{\bm\nu}
\newcommand{\ep}{\epsilon}
\newcommand{\ovl}{\overline}
\newcommand{\de}{\delta}
\begin{document}

\title[Simultaneous non-vanishing of  $GL(2) \times GL(3) $ and $GL(3)$ $L$-functions]{Simultaneous non-vanishing of central values of $GL(2) \times GL(3) $ and $GL(3)$ $L$-functions}
\author{Gopal Maiti and Kummari Mallesham}
\address{ Kummari Mallesham \newline {\em Stat-Math Unit, Indian Statistical Institute, 203 B.T. Road, Kolkata 700108, India; \newline  Email: iitm.mallesham@gmail.com
} }

\address{ Gopal Maiti \newline {\em Stat-Math Unit, Indian Statistical Institute, 203 B.T. Road, Kolkata 700108, India; \newline  Email: g.gopaltamluk@gmail.com
} }

\maketitle

\begin{abstract}
We study simultaneous non-vanishing of $L(\tfrac{1}{2},\di)$ and $L(\tfrac{1}{2},g\otimes \di)$, when $\di$ runs over an orthogonal basis of the space of Hecke-Maass cusp forms for $SL(3,\mathbb{Z})$ and $g$ is a fixed $SL(2,\mathbb{Z})$ Hecke cusp form of weight $k\equiv 0 \pmod 4$.
\end{abstract}
\blfootnote{
	{\it key words and phrases}: Rankin-Selberg $L$-functions, $GL(3)$ $L$-function, Kuznetsov trace formula,\\ Simultaneous non-vanishing.}

\section{Introduction}
Like the Birch and Swinnerton-Dyer conjecture which relates the order of vanishing of the Hasse-Weil $L$-function at the central point to the rank of an elliptic curve, the vanishing or non-vanishing of an automorphic $L$-function at the special points are related to several deep problems with great significance in number theory. Therefore, it is a profoundly interesting question to understand whether product of two or more $L$-functions are simultaneous non-vanishing at the central point. This type of question has been studied by many authors (see for example \cite{IS}, \cite{MV}, \cite{KMV}, \cite{MS}, \cite{Liu2}, \cite{HLX}, \cite{Go}). In $2014$, Das \& Khan \cite{DK} proved that $GL(2)\times GL(1)$ and $GL(1)$ $L$-functions are simultaneous non-vanishing. Ramakrishnan \& Rogawski \cite{RR} in $2005$, showed a simultaneous non-vanishing result for $GL(2)\times GL(1)$ and $GL(2)$ $L$-functions. Similar type of non-vanishing  results for $GL(2)\times GL(2)$ and $G(2)$ $L$-functions were proved by Xu \cite{Xu} and Liu \cite{Liu1} for Maass form and holomorphic Hecke cusp form with weight aspect respectively. 
Non-vanishing problem for $GL(3)\times GL(2)$ and $GL(2)$ $L$-functions was first studied by Li \cite{Li} in $2009$. More precisely, let $f$ be a fixed Hecke-Maass cusp form for $SL(3)$. Li Proved that there are infinitely many $SL(2)$ Hecke-maass cusp forms $u_j$ such that $L(\frac{1}{2},f\times u_j )\,L(\frac{1}{2}, u_j)\neq 0$. In the level aspect for holomorphic Hecke cusp form case similar result was obtained by Khan in \cite{Kh}.

\vspace{1 mm}
 In this paper, we consider the first moment of the product of $GL(2)\times GL(3)$ and $GL(3)$ $L$-functions. More precisely, we fix an Hecke cusp form $g$ for $SL(2,\mathbb{Z})$. 
 We study first moment of $L(s,\di\otimes g) L(s,\di)$ at the central point as $\di$ runs over an orthogonal basis of the space of Hecke-Maass cusp forms for $SL(3,\mathbb{Z})$.  

\noindent
 Let $\bmu_\di=(\mu_{1},\mu_{2},\mu_{3})$ be the Langlands parameter and $\bnu_{\di}=(\nu_{1},\nu_{3},\nu_{3})$ be the spectral parameter of a Hecke-Maass cusp form $\di$ for $SL(3,\mathbb{Z})$. As in Blomer-Buttcane \cite{BB}, we consider the generic case in short interval. Let $\bmu_0=(\mu_{0,1},\mu_{0,2},\mu_{0,3})$ be a fixed point in $\Lambda^{\prime}_{1/2}$ (see \eqref{Lamda}). So $\bnu_{0}=(\nu_{0,1},\nu_{0,2},\nu_{0,3})$ satisfies the relations \eqref{mu nu rel1} and \eqref{mu nu rel2}. We consider the case 
 \[ 
 |\mu_{0,j}|\asymp|\nu_{0,j}|\asymp\|\bmu_0\|\asymp\|\bnu_0\|:=T,\;\; 1\le j\le 3. 
 \]
Let us denote $R=T^{\theta}$ for any fixed $\theta$ in $(0,1)$. We choose the test function $h(\bmu)$ so that it has the localizing effect at a ball of radius $R$ about $w(\bmu_0)$, where $w$ are elements in the Weyl group $\mathfrak{W}$ of $GL(3,\mathbb{R})$.  It is defined by 
\[
h(\bmu):=P(\bmu)^{2}\left( \sum_{w\in \mathcal{W}}\psi\left(\frac{w(\bmu)-\bmu_0}{R} \right) \right)^{2} ,
\]
where $\psi(\bmu)=\exp\left(-(\mu_1^{2}+\mu_2^{2}+\mu_3^{2})\right)$ and 
\[
    P(\bmu)=\prod_{1\le n\le A_0 }\prod_{j=1}^{3}\frac{\left(\nu_j -\tfrac{1}{3}(1+2n)\right) \left(\nu_j +\tfrac{1}{3}(1+2n)\right)}{|\nu_{0,k}|^2}
\] for some fixed large $A_0 >0$. Here 
\[
\mathfrak{W}:=\left\{I, w_2 =\left(\begin{smallmatrix}
1 & &\\
  & & 1\\
  & 1 &
\end{smallmatrix}\right), w_3 =\left(\begin{smallmatrix}
& 1 &\\
1 & &\\
& & 1
\end{smallmatrix}\right), w_4 =\left(\begin{smallmatrix}
& 1 &\\
& & 1\\
1 & &
\end{smallmatrix}\right), w_5  =\left(\begin{smallmatrix}
& & 1\\
1 & &\\
& 1 &
\end{smallmatrix}\right), w_6 =\left(\begin{smallmatrix}
& & 1\\
& 1 &\\
1 & &
\end{smallmatrix}\right)    \right\}
\] is the Weyl group for $SL(3,\mathbb{R})$. 
 Let $ \rm  d_{\text{spec}}\bmu =\text{spec}(\bmu)d\bmu$ with $$\rm  \text{spec}(\bmu)=\displaystyle\prod_{j=1}^{3}\left(3\nu_j \tan\left(\frac{3\pi}{2}\nu_j\right)\right)\quad \text{and} \quad d\bmu=d\mu_1\,d\mu_2=d\mu_2\,d\mu_3=d\mu_3\,d\mu_1.$$ Let us define $$\mathcal{N}_{\di}=\|\di\|^{2}\displaystyle\prod_{j=1}^{3}\cos\left(\frac{3\pi}{2}\nu_j\right)$$  to be the normalizing factor.
 \noindent
 Now we state the main theorem of this article. 
\begin{theorem}\label{Theo}
	Let $g$ be a Hecke cusp form for $SL(2,\mathbb{Z})$ of weight $k\equiv 0 \pmod 4$. Let $\{\di\}$ be a basis of the space of Hecke-Maass cusp forms for $SL(3,\mathbb{Z})$. Then we have
	\begin{align*}
	\sum_{\di}\frac{h(\boldsymbol{\mu}_{\di})}{\mathcal{N}_{\di}} L(\tfrac{1}{2},g\otimes \di) L(\tfrac{1}{2}, \di)=\frac{1}{192\,\pi^{5}}\iint_{\Re(\bm\mu)=0}M(\bm\mu,k) h(\bm{\mu})\mathrm{spec}(\bm\mu)\rm d\bm\mu\\
	 +O\left(T^{\frac{17}{6}+\varepsilon}R^{2}\right),\hspace{2cm}
	\end{align*}
where 
\begin{align*}
M(\bm\mu,k)=\zeta(\tfrac{3}{2})+L(1,g)\prod_{j=1}^{3}\frac{\Gamma(\tfrac{1}{4}+\tfrac{\mu_{j}}{2})}{\Gamma(\tfrac{1}{4}-\tfrac{\mu_{j}}{2})} + L(1,g)\prod_{j=1}^{3}\frac{\Gamma(\tfrac{k}{2}+\mu_{j})}{\Gamma(\tfrac{k}{2}-\mu_{j})}\\
+ \zeta(\tfrac{3}{2})\prod_{j=1}^{3} \frac{\Gamma(\tfrac{k}{2}+\mu_{j})\Gamma(\tfrac{1}{4}+\tfrac{\mu_{j}}{2})}{\Gamma(\tfrac{k}{2}-\mu_{j})\Gamma(\tfrac{1}{4}-\tfrac{\mu_{j}}{2})}.
\end{align*}
\end{theorem}
 
\noindent
Note that, $\displaystyle\iint_{\Re(\bm\mu)=0}M(\bm\mu,k) h(\bm{\mu})\mathrm{spec}(\bm\mu)d\bm\mu\asymp T^{3}R^{2}$.
\begin{remark}
We have used bounds $zeta(1/2 +it)\ll t^{1/6}$ and $L(1/2 +it,g \otimes f) \ll t$ in the estimation of the error terms (see section\ref{S3.16}). Note that by using best know bounds for $\zeta(1/2 +it)$ and $L(1/2 +it,g \otimes f)$ one can get slight improvement in the error term. Since we are only focusing on simultaneous non-vanishing results we have not incorporate such a small improvement in the error term.
\end{remark}
\noindent
 As a corollary of Theorem \ref{Theo}, we have the following result.	  
\begin{corollary} \label{Cor1}
	Let $g$ be a Hecke cusp form of weight $k\equiv 0\pmod 4$ for $SL(2,\mathbb{Z})$. Then there exist infinitely many Hecke-Maass cusp forms $\di$ for $SL(3,\mathbb{Z})$ such that $L(\tfrac{1}{2}, \di)  L(\tfrac{1}{2},g\otimes \di)\neq 0$.
\end{corollary}
 
 \begin{remark}
 	If $g$ is a Hecke cusp form for $SL(2,\mathbb{Z})$ with $k\equiv 2 \mod4$, then $L(\tfrac{1}{2},g\ot\di)=0$. In that case one can consider the following sum \[\sum_{\di}\frac{h(\boldsymbol{\mu}_{\di})}{\mathcal{N}_{\di}} L^{\prime}(\tfrac{1}{2},g\otimes \di) L(\tfrac{1}{2}, \di) \]
 	and get similar results.
 \end{remark}

\vspace{2mm}
\section{Preliminaries} In this section we review some definitions, essential facts and tools that will be used in later development.

\subsection{Automorphic forms for $SL(3,\mathbb{Z})$ and their $L$-functions.} 
Let $$\mathbb{H}_{3}=  GL(3,\mathbb{R})/ O(3,\mathbb{R}) \mathbb{R}^{*}$$ be the generalized upper half plane. For $0\le c\le \infty$, let
\begin{equation}\label{Lamda}
\Lambda_{c}'=\left\{\bmu= (\mu_1,\mu_2,\mu_3)\in\mathbb{C}^3,\quad
\begin{aligned}
&|{\Re}(\mu_j)|\le c,\quad  \mu_1+\mu_2+\mu_3=0,\\
&\{-\mu_1,-\mu_2,-\mu_3\}=\{\overline\mu_1,\overline\mu_2,\overline\mu_3\}
\end{aligned}
\right\}.
\end{equation}
Consider $\bmu$ to be the Langlands parameter of a Hecke-Maass form $\di$ in $L^{2}\left(SL(3,\mathbb{Z})\backslash \mathbb{H}_{3}\right)$. Let us define
\begin{equation}\label{mu nu rel1}
\nu_1=\frac{1}{3}(\mu_1 -\mu_2),\quad \nu_2=\frac{1}{3}(\mu_2 -\mu_3),\quad\nu_3 =-\nu_1 -\nu_2=\frac{1}{3}(\mu_3 -\mu_1)
\end{equation}
where ${\bm \nu}=(\nu_1 ,\nu_2 ,\nu_3)$ is known as the spectral parameter of $\di$. So
\begin{align}\label{mu nu rel2}
\mu_1= 2\nu_1+ \nu_2 , \quad\mu_2= \nu_2- \nu_1,\quad \mu_3=-\nu_1-2\nu_2.
\end{align}

\noindent
Let $A_{\di}(m_{1},m_{2})$ be the normalised Fourier coefficients of a $GL(3)$ Hecke Maass cusp form $\di$ with Langlands parameters $\bm\mu_\di= (\mu_{1},\mu_{2}, \mu_{3})$. The stander $L$-function associated to $\di$ is given by
$$L(s,\di) = \sum_{n \geq 1} \frac{A_{\di}(1,n)}{n^{s}}\quad \text{for}\;\, \Re{(s)} >1.$$
The dual form of $\di$ is denoted by $\widetilde{\di}$ with the Langlands parameter $\bmu_{\tilde\di}=(-\mu_{1},-\mu_{2}, -\mu_{3})$ and the coefficients $A_{\di}(n,1)=\overline{ A_{\di}(1,n)}=A_{\tilde{\di}}(1,n)$. Let us define $$\Lambda(s,\di):= \gamma(s,\di) \, L(s, \di), $$
where $\gamma(s, \di )=\displaystyle \prod_{j=1}^{3}  \Gamma_{\mathbb{R}}\left(s-\mu_{j} \right)$ and $\Gamma_{\mathbb{R}}(s)=\pi^{-\tfrac{s}{2}}\Gamma(\tfrac{s}{2})$. 
$\Lambda(s,\di)$ is called the completed $L$-function,
which is an entire function and satisfies the functional equation
$$\Lambda(s,\di) = \Lambda(1-s, \tilde{\di}).$$

\subsection{The maximal Eisenstein series.}\label{Max Eisen}
Let $u\in \mathbb{C}$ have sufficiently large real part. Let $f$ be a Hecke-Maass cusp form for $SL(2,\mathbb{Z})$ with the spectral parameter $i t_{f}$, Hecke eigenvalues $\lambda_{f}(m)$ and $\|f\|=1$. The maximal Eisenstein series and it's Hecke eigenvalue at $(m,n)$ are denoted by $E_{u,f}^{\max}(z)$ and $B_{u,f}^{\max}(m,n)$, respectively. The Hecke eigenvalue $B_{u,f}^{\max}(1,m)$  is defined by (see Goldfeld \cite{Gol}) 
\begin{align*}
B_{u,f}^{\max}(1,m)=\sum_{d_1 d_2 =m} \lambda_{f}(d_1) d_{1}^{-u} d_{2}^{-2u}
\end{align*}
and satisfies the following Hecke relations
\begin{align*}
B_{u,f}^{\max}(m,1)=\overline{B_{u,f}^{\max}(1,m)},\quad
B_{u,f}^{\max}(m,n)=\sum_{d\mid (m,n)}\mu(d) B_{u,f}^{\max}(\tfrac{m}{d},1) B_{u,f}^{\max}(1,\tfrac{n}{d}).
\end{align*}
The $L$-function associated to $E_{u,f}^{\max}(z)$ is given by
\begin{align}\label{Emax}
L(s,E_{u,f}^{\max})=\sum_{m\ge 1}\frac{B_{u,f}^{\max}(1,m)}{m^s}=\zeta(s-2u) L(s+u,f),
\end{align}
for sufficiently large $\Re(s)$. 
It satisfies the functional equation 
\begin{align*}
\Lambda(s,E_{u,f}^{\max})=\prod_{j=1}^{3}\Gamma_{\mathbb{R}}(s+\mu_{j}^{\prime})L(s,E_{u,f}^{\max})=\Lambda(1-s,E_{-u,f}^{\max}),
\end{align*}
where $\mu_{1}^{\prime}=u+it_{f}$, $\mu_{2}^{\prime}=u-it_{f}$ and $\mu_{3}^{\prime}=-2u$.
 \noindent
The normalized factor for the maximal Eisenstein series is define by 
\[ 
    \mathcal{N}_{u,f}^{\max}:= 8L(1,{\rm Ad}^{2}f) |L(1+3u,f)|^{2}.  
\]

\subsection{The minimal Eisenstein series.}\label{Min Eisen}
Let $\nu_{1}, \nu_2 \in \mathbb{C}$ and $(\mu_1 ,\mu_2 ,\mu_3)$ be the Langlands parameter given by \eqref{mu nu rel2}. We denote the minimal Eisenstein series by $E_{\nu_1 ,\nu_2}^{\min}(z)$. The Hecke eigenvalue $B_{\nu_1 ,\nu_2}^{\min}(m,n)$ of $E_{\nu_1 ,\nu_2}^{\min}(z)$  at $(m,n)$ is defined by (see Goldfeld \cite{Gol})
\[
B_{\bmu}^{\rm min}(1,m)=B_{\nu_1 ,\nu_2}^{\min}(1,m):=\sum_{d_{1} d_{2} d_{3}=n}d_{1}^{-\mu_1}d_{2}^{-\mu_2}d_{3}^{-\mu_3}
\]
and satisfies the following Hecke relations
\begin{align*}
B_{\nu_1 ,\nu_2}^{\min}(m,1)=\overline{B_{\nu_1 ,\nu_2}^{\min}(1,m)},\quad
B_{\nu_1 ,\nu_2}^{\min}(m,n)=\sum_{d\mid (m,n)}\mu(d) B_{\nu_1 ,\nu_2}^{\min}(\tfrac{m}{d},1) B_{\nu_1 ,\nu_2}^{\min}(1,\tfrac{n}{d}).
\end{align*}
The $L$-function associated to $E_{\nu_1 ,\nu_2}^{\max}(z)$ is given by
\begin{align}\label{Emin}
L(s,E_{\nu_1 ,\nu_2}^{\min})=\sum_{m\ge 1}\frac{B_{\nu_1 ,\nu_2}^{\min}(1,m)}{m^s}=\zeta(s+\mu_1 ) \zeta(s+\mu_2 )\zeta(s+\mu_3 ),
\end{align}
for $\Re(s)>1$. The normalized factor for the minimal Eisenstein series is define by 
\[ 
\mathcal{N}_{\bmu}^{\rm min}=\mathcal{N}_{\nu_1 , \nu_2}^{\rm min}:= \frac{1}{16}\prod_{j=1}^{3} |\zeta(1+3\nu_j)|^{2}.  
\]

\subsection{The Rankin-Selberg $L$-function on GL(2) $\times$ GL(3).}
Let $\di$ be a $GL(3)$ Hecke Maass cusp with Langlands parameters $\mu= (\mu_{1},\mu_{2}, \mu_{3})$. Let $g$ be a Hecke cusp form for $SL(2,\mathbb{Z})$ of weight $k$ and $\lambda_{g}(n)$ be the $n$-th Hecke eigenvalue. The Rankin-Selberg  $L$-function of $g$ and $\di$ is defined by
\begin{equation*}
L\left(s, g\otimes \di\right) = \mathop{\sum \sum}_{m,n\geq 1}  \frac{\lambda_{g}(n) \overline{A_{\digamma}(m,n)}}{\left(nm^2\right)^{s}} \quad, \, \Re{(s)} >1.
\end{equation*}
 It is entire and satisfies the functional equation
 $$\Lambda(s, g\otimes \di)= i^{3k} \, \Lambda(1-s,g\otimes\tilde{\di}) ,$$
 where 
 $$\Lambda(s, g\otimes \di) = \gamma(s,g\otimes \di) \, L(s, g\otimes \di)$$
and
 $$\gamma(s, g\otimes \di)= \prod_{j=1}^{3} \Gamma_{\mathbb{R}} \left(s+ \frac{k-1}{2}-\mu_{j} \right)  \Gamma_{\mathbb{R}} \left(s+ \frac{k+1}{2}-\mu_{j} \right).$$
 
\noindent
Let $E_{u,f}^{\max}$ be the maximal Eisenstein series as in \S \ref{Max Eisen}. The Rankin-Selberg $L$-function $L(s, g\otimes E_{u,f}^{\max} )$ is defined by 
\begin{align}\label{gEmax} 
L(s, g\otimes E_{u,f}^{\max})&=\mathop{\sum \sum}_{m,n\geq 1}  \frac{\lambda_{g}(n) \overline{B_{u,f}^{\max}(m,n)}}{\left(nm^2\right)^{s}}\nonumber\\
&= L(s+2u,g)L(s-u,g\otimes f),
 \end{align}
for sufficiently large $\Re(s)$.

\noindent
Let $E_{\nu_{1},\nu_{2}}^{\min}$  be the minimal Eisenstein series as in \S \ref{Min Eisen}. The Rankin-Selberg $L$-function $L(s, g\otimes E_{\nu_{1},\nu_{2}}^{\min} )$ is defined by
 \begin{align}\label{gEmin} 
 L(s, g\otimes E_{\nu_{1},\nu_{2}}^{\min})&=\mathop{\sum \sum}_{m,n\geq 1}  \frac{\lambda_{g}(n) \overline{B_{\nu_{1},\nu_{2}}^{\min}(m,n)}}{\left(nm^2\right)^{s}}\nonumber\\
 &= L(s-\mu_1 ,g)L(s-\mu_2 , g)L(s-\mu_3 , g)
 \end{align}
 for  $\Re(s)>1$.
 
\vspace{2mm}  
\noindent 
 From the above functional equation we deduce the approximate functional equation for $L(s,\di)$ and $L(s,g\ot \di)$ at the central point $s=\frac{1}{2}$ (See \S $5.2$ of \cite{IK}).
 Let $G(s)=e^{s^2}$. We define $$V_{\di}(y):= \frac{1}{2 \pi i} \int_{(3)} y^{-u} \frac{\gamma(\tfrac{1}{2}+u, \di)}{\gamma(\tfrac{1}{2}, \di )} G(u) \frac{{\rm d}u}{u},$$
 $$ \tilde{V}_{\di}(y): = \frac{1}{2 \pi i} \int_{(3)} y^{-u} \frac{\gamma(\tfrac{1}{2}+u, \tilde{\di} )}{\gamma(\tfrac{1}{2}, \di)} G(u) \frac{{\rm d}u}{u},$$
 and
 $$W_{\di}(y): = \frac{1}{2 \pi i} \int_{(3)} y^{-u} \frac{\gamma(\tfrac{1}{2}+u, g\otimes \di)}{\gamma(\tfrac{1}{2}, g\otimes \di)} G(u) \frac{{\rm d}u}{u},$$
 $$
\tilde{W}_{\di}(y): = \frac{1}{2 \pi i} \int_{(3)} y^{-u} \frac{\gamma(\tfrac{1}{2}+u, g\otimes\tilde{\di})}{\gamma(\tfrac{1}{2}, g\otimes \di)} G(u) \frac{{\rm d}u}{u}.$$
 
 \begin{lemma}\label{le1}
 	 We have 
 	$$L(\tfrac{1}{2},\di) = \mathop{\sum }_{n \geq 1} \frac{A_{\di}(1,l)}{l^{1/2}} V_{\di}(l) + \mathop{\sum }_{l \geq 1} \frac{\overline{A_{\di}(1,l)}}{l^{1/2}} \tilde{V}_{\di}(l).$$	
 Moreover, for $\bmu=(\mu_{1}, \mu_{2}, \mu_{3})$ with $\mu_j \asymp T$, one has
 \[
 y^{j_1}\frac{\partial^{j_1}}{\partial y^{j_1}}V_{\di}(y)\ll \left(\frac{y}{T^3}\right)^{-A},\quad y^{j_1}\frac{\partial^{j_1}}{\partial y^{j_1}}\tilde{V}_{\di}(y)\ll \left(\frac{y}{T^3}\right)^{-A} 
 \]	
 for any $A>0$ and any $j_1 \in \mathbb{N}\cup\{0\}$. Also, for $y\gg T^3$ 
 \[
 V_{\di}(y)=1+O_{A}\left(\frac{T^3}{y}\right)^{-A},\quad  \tilde{V}_{\di}(y)=\prod_{j=1}^{3}\frac{\Gamma\left(\tfrac{1}{4}+\tfrac{\mu_j}{2}\right)}{\Gamma\left(\tfrac{1}{4}-\tfrac{\mu_j}{2}\right)}+O_{A}\left(\frac{T^3}{y}\right)^{-A}
 \]	for any $A>0$.
 \end{lemma}
 
 \begin{lemma}\label{le2} We have 
 	$$L(\tfrac{1}{2}, g\otimes \di) = \mathop{\sum \sum}_{m,n \geq 1} \frac{\lambda_{g}(n) A_{\di}(m,n)}{\left(n m^2\right)^{1/2}} \, W_{\di}(nm^2) + i^{3k} \mathop{\sum \sum}_{m,n \geq 1} \frac{\lambda_{g}(n) \overline{A_{\di}(m,n)}}{\left(n m^2\right)^{1/2}} \tilde{W}_{\di}(nm^2),$$ 	
 	 Moreover, for $\bmu=(\mu_{1}, \mu_{2}, \mu_{3})$ with $\mu_j \asymp T$, one has
 	\[
 	y^{j_1}\frac{\partial^{j_1}}{\partial y^{j_1}}W_{\di}(y)\ll_{k} \left(\frac{y}{T^3}\right)^{-A},\quad y^{j_1}\frac{\partial^{j_1}}{\partial y^{j_1}}\tilde{W}_{\di}(y)\ll_{k} \left(\frac{y}{T^3}\right)^{-A} 
 	\]	
 	for any $A>0$ and any $j_1 \in \mathbb{N}\cup\{0\}$. Also, for $y\gg T^3$ 
 	\[
 	W_{\di}(y)=1+O_{A, k}\left(\frac{T^3}{y}\right)^{-A},\quad  \tilde{W}_{\di}(y)=\prod_{j=1}^{3}\frac{\Gamma\left(\tfrac{k}{2}+\mu_j\right)}{\Gamma\left(\tfrac{k}{2}-\mu_j\right)}+O_{A, k}\left(\frac{T^3}{y}\right)^{-A}
 	\]	for any $0<A<\frac{k-1}{2}$.
 	
 \end{lemma}

\subsection{The Kloostermam sums.}
For $n_1 , n_2 , m_1 , m_2 , D_1 , D_2 \in \mathbb{N}$, we define the following Kloosterman sums.
\[
\tilde{S}(n_1 , n_2 , m_1 ; D_1 , D_2 ):=\mathop{\sum \sum}_{\substack{C_1({\rm mod}\, D_1 ),\, C_2({\rm mod}\, D_2 ) \\ (C_1 ,D_1)=(C_2 ,D_{2}/D_{1})=1 }} e\left(n_2 \frac{\bar{C_{1}} C_2}{D_1} + m_1 \frac{\bar{C_{2}}}{D_{2}/D_{1}} + n_1 \frac{C_1}{D_1} \right)
\]
for $D_1 \mid D_2$, and  

\begin{align*}
&S(n_1 ,m_2 , m_1 , n_2 ; D_1 , D_2 ) \\
&:=\mathop{\sum \sum \sum \sum}_{\substack{B_1 , C_1({\rm mod}\, D_1 );\, B_2 , C_2({\rm mod}\, D_2 )\\ D_1 C_2 + B_1 B_2 + D_2 C_1 \equiv 0\, ({\rm mod}\, D_1 D_2 ) \\ (B_j  ,C_j ,D_j)=1 }} e\left( \frac{n_1 B_1 + m_1 (Y_1 D_2 - Z_1 B_2)}{D_1} + \frac{m_2 B_2 + n_2(Y_2 D_1 - Z_2 B_1)}{D_2} \right),
\end{align*}
where $  B_j Y_j + C_j  Z_j\equiv 1\, ({\rm mod}\, D_j )$ for $j=1,2$.

 \subsection{Integral Kernels.} 
Following [\cite{Bt}, Theorem $2$ \& $3$], we define the integral kernel in terms of Mellin-Barnes representations. For $s\in \mathbb{C}$, $\bmu=(\mu_1 , \mu_2 , \mu_3)$ define the meromorphic function 
\[  
\tilde{G}^{\pm}(s,\bmu):=\frac{\pi^{-3s}}{12288\pi^{7/2}}\left(\prod_{j=1}^{3}\frac{\Gamma(\tfrac{1}{2}(s-\mu_j))}{\Gamma(\frac{1}{2}(1-s+\mu_j))}\pm i \prod_{j=1}^{3}\frac{\Gamma(\frac{1}{2}(1+s-\mu_j))}{\Gamma(\frac{1}{2}(2-s+\mu_j))} \right),
\]
and for ${\bm s}=(s_1 ,s_2)\in \mathbb{C}^2 $, $\bmu=(\mu_1 , \mu_2 , \mu_3)$ define the meromorphic function 
\[
G(\bm s, \bmu):=\frac{1}{\Gamma(s_1 +s_2)} \prod_{j=1}^{3} \Gamma(s_1 -\mu_j )\Gamma(s_2 +\mu_j ).
\]
We also define the following  trigonometric functions 
\begin{align*}
S^{++}(\bm s; \bmu)&:=\frac{1}{24\pi^{2}}\prod_{j=1}^{3}\cos\left(\tfrac{3}{2}\pi \nu_j\right),\\
S^{+-}(\bm s; \bmu)&:=-\frac{1}{32\pi^{2}}\frac{\cos\left(\tfrac{3}{2}\pi \nu_2\right)\sin\left(\pi(s_1 -\mu_1)\right)\sin\left(\pi(s_2 +\mu_2)\right)\sin\left(\pi(s_2 +\mu_3)\right)}{\sin\left(\tfrac{3}{2}\pi \nu_1\right)\sin\left(\tfrac{3}{2}\pi \nu_3\right)\sin\left(\pi(s_1 + s_2)\right)},\\
S^{-+}(\bm s; \bmu)&:=-\frac{1}{32\pi^{2}}\frac{\cos\left(\tfrac{3}{2}\pi \nu_1\right)\sin\left(\pi(s_1 -\mu_1)\right)\sin\left(\pi(s_1 -\mu_2)\right)\sin\left(\pi(s_2 +\mu_3)\right)}{\sin\left(\tfrac{3}{2}\pi \nu_2\right)\sin\left(\tfrac{3}{2}\pi \nu_3\right)\sin\left(\pi(s_1 + s_2)\right)},\\
S^{--}(\bm s; \bmu)&:=\frac{1}{32\pi^{2}}\frac{\cos\left(\tfrac{3}{2}\pi \nu_3\right)\sin\left(\pi(s_1 -\mu_2)\right)\sin\left(\pi(s_2 +\mu_2)\right)}{\sin\left(\tfrac{3}{2}\pi \nu_2\right)\sin\left(\tfrac{3}{2}\pi \nu_1\right)}.
\end{align*}
For $y\in \mathbb{R}^{*}$ with ${\rm sgn} (y)=\ep $, let 
\[
K_{w_4}(y;\bmu):=\int_{-i \infty}^{+i \infty} |y|^{-s}\tilde{G}^{\ep}(s,\bmu)\frac{{\rm d}s}{2\pi i}.
\]
For ${\bm y}=(y_1 ,y_2)\in{(\mathbb{R}^{*})}^{2}$ with ${\rm sgn} (y_1)=\ep_1$, ${\rm sgn} (y_2)=\ep_2$, let 
\[
K_{w_6}^{\ep_1 ,\ep_2}(\bm y;\bmu):=\int_{-i \infty}^{+i \infty}\int_{-i \infty}^{+i \infty} |4\pi^2 y_1|^{-s_1} |4\pi^{2} y_2|^{-s_2} G(\bm s, \bmu)S^{\ep_1 ,\ep_2}(\bm s;\bmu)\frac{{\rm d}s_1 {\rm d}s_2}{(2\pi i)^2}.
\]

\subsection{The Kuznetsov formula} Define the spectral measure on the hyperplane $\mu_{1}+\mu_{2}+\mu_{3}=0$ by $\rm d_{ spec}\bmu ={spec}(\bmu) d\bmu$, where

$\rm\text{spec}(\bmu)=\displaystyle\prod_{j=1}^{3}\left(3\nu_j \tan\left(\frac{3\pi}{2}\nu_j\right)\right)\quad \text{and} \quad d\bmu=d\mu_1\,d\mu_2=d\mu_2\,d\mu_3=d\mu_3\,d\mu_1.$

\noindent
 Now we state the Kuznetsov trace formula in the version of Buttcane [\cite{Bt}, Theorem $2, 3, 4$].

\begin{lemma}\label{Kuz}
	Let $n_1 , n_2 , m_1 , m_2 \in \mathbb{N}$ and $h$ be a holomorphic function on 
	$$\Lambda_{\frac{1}{2}+\delta}=\left\{\bmu=(\mu_1 , \mu_2 , \mu_3)\in\mathbb{C}^{3}, \mu_1 +\mu_2 + \mu_3 =0, \Re{(\mu_j)}\le \frac{1}{2}+\delta  \right\}$$
	for some $\delta >0$, symmetric under the Weyl group $\mathcal{W}$, of rapid decay when $|\Im(\mu_j)|\to \infty$ and satisfies 
	$$h(3\nu_1 \pm 1, 3\nu_2 \pm 1, 3\nu_3 \pm 1)=0.$$
	Then we have $$\mathcal{C}+\mathcal{E}_{\min}+\mathcal{E}_{\max}=\Delta+\Sigma_4 +\Sigma_{5} +\Sigma_{6},$$
	where
	\begin{align*}
	\mathcal{C}&:=\sum_{\di} \frac{h(\bm\mu_{\di})}{\mathcal{N}_{\di}} \ovl{A_{\di}(m_1 ,m_2 )} A_{\di}(n_1 ,n_2 ),\\
	\mathcal{E}_{\rm min}&:=\frac{1}{24 (2 \pi i)^2} {\iint}_{\Re(\bmu)=0} \frac{h(\bmu)}{\mathcal{N}_{\bmu}^{\rm min}} \ovl{B_{\bmu}^{\min}(m_1,m_2)} B_{\bmu}^{\min}(n_1,n_2)\, \rm d\mu_1 d\mu_2,\\
	\mathcal{E}_{\rm max}&:= \sum_{f}\frac{1}{2 \pi i} \int_{\Re(u)=0}\frac{h(u+it_f , u-it_f , -2u )}{\mathcal{N}_{u,f}^{\rm max}}\ovl{B_{u,f}^{\rm max}(m_1 ,m_2)} B_{u,f}^{\rm max}(n_1 ,n_2) {\rm d}u,
	\end{align*}
	and 
	\begin{align*}
	\Delta&:=\delta_{m_1 , n_1} \delta_{m_2 , n_2} \frac{1}{192 \pi^{5}}{\iint }_{\Re(\bmu)=0} h(\bmu)\, \rm d_{spec}\bmu,\\
	\Sigma_4 &:=\sum_{\ep=\pm 1}\sum_{\substack{D_2 \mid D_1 \\ m_2 D_1 =n_1 {D_2}^2}}\frac{\tilde{S}(-\ep n_2 , m_2 , m_1 ; D_2 , D_1)}{D_1 D_2} \Phi_{w_4}\left(\frac{\ep m_1 m_2 n_2}{D_1 D_2}\right),\\
	\Sigma_5 &:=\sum_{\ep=\pm 1}\sum_{\substack{D_1 \mid D_2 \\ m_1 D_2 =n_2 {D_1}^2}}\frac{\tilde{S}(\ep n_1 , m_1 , m_2 ; D_1 , D_2)}{D_1 D_2} \Phi_{w_5}\left(\frac{\ep m_1 m_2 n_1}{D_1 D_2}\right),\\
	\Sigma_6 &:=\sum_{\ep_1 ,\ep_2 =\pm 1}\sum_{D_1 ,D_2}\frac{S(\ep_2 n_2 ,\ep_1 n_1 , m_1 , m_2 ; D_1 , D_2)}{D_1 D_2} \Phi_{w_6}\left(-\frac{\ep_2 m_1 n_2 D_2}{{D_1}^{2}},-\frac{\ep_1 m_2 n_1 D_1}{{D_2}^{2}} \right) 
	\end{align*}
	with 
	\begin{align}\label{Phi}
	\Phi_{w_4}(y)&:=\iint_{\Re(\bmu)=0} h(\bmu)K_{w_4}(y;\bmu)\rm d_{spec}\bmu, \nonumber\\ 
	\Phi_{w_5}(y)&:=\iint_{\Re(\bmu)=0} h(\bmu)K_{w_4}(- y;-\bmu)\rm d_{spec}\bmu,\\ 
	\Phi_{w_6}(\bm y)&:=\iint_{\Re(\bmu)=0} h(\bmu)K_{w_6}^{\rm sgn(y_1), sgn(y_2)}(\bm y;\bmu)\rm d_{spce}\bmu.\nonumber
	\end{align}
\end{lemma}

Here we quote Lemma $8$ and Lemma $9$ of \cite{BB} which are used in truncating summation in geometric terms after the application of the Kuznetsov formula.

\begin{lemma}\label{le4}
	 Let $0<|y|\le T^{3-\ve}$. Then for any constant $B>0$, we have 
	\[ 
	  \Phi_{w_4}(y)\ll_{\ve, B} T^{-B}.
	\]
If $|y|>T^{3-\ve}$, then 
\[ 
|y|^{j}\frac{\rm d^{j}\,}{\;{\rm d}y^{j}}\Phi_{w_4}(y)\ll_{\ve, j} T^{3+\ve}R^{2}\left(T+|y|^{1/3}\right)^{j}
\]
for any $j\in \mathbb{N}\cup\{0\}$.

\end{lemma}

\begin{lemma}\label{le5}
	 Let $\Upsilon:=\min\left\{|y_1 |^{1/3}|y_2 |^{1/6} , |y_2 |^{1/3}|y_1 |^{1/6} \right\}$. If $\Upsilon\ll T^{1-\ve}$, then for any constant $B>0$
	\[
	\Phi_{w_6}(y_1 , y_2 )\ll_{\ve, B} T^{-B}.
	\]
If $\Upsilon\gg T^{1-\ve}$, then we have
\begin{align*}
&|y_1|^{j_1}|y_2|^{j_2}\frac{\partial^{j_1} }{\partial y_{1}^{j_1}} \frac{\partial^{j_2} }{\partial y_{2}^{j_2}}\Phi_{w_6}(y_1 , y_2)\\
& \ll_{\ve, j_1 ,j_2} 
T^{3}R^{2}\left( T + |y_1|^{1/2} + |y_1 |^{1/3} |y_2 |^{1/6} \right)^{j_1}\left( T + |y_2|^{1/2} + |y_2 |^{1/3} |y_1 |^{1/6} \right)^{j_2}
\end{align*}
for any $j_1 ,j_2 \in \mathbb{N}\cup\{0\}$.	
\end{lemma}

\section{proof of theorem \ref{Theo}}
Let us define
$$\mathcal{S}(T) := \sum_{\di} \frac{h(\bm\mu_{\di})}{\mathcal{N}_{\di}} \, L(1/2, g\otimes \di) \, L(1/2, \di).$$
By using the approximate functional equations of $L(1/2, g\otimes \di  )$ and $L(1/2,\di)$ (Lemma \ref{le2} and Lemma \ref{le1} ), we get
 \begin{align*}
 \mathcal{S}(T) &= \mathop{\sum \sum}_{m,n \geq 1} \sum_{l \geq 1} \frac{\lambda_{g}(n)}{m \left( nl\right)^{1/2}}  \sum_{\di} \frac{h_{1}(\bm\mu_{\di})}{\mathcal{N}_{\di}} A_{\di}(m,n) \, \overline{A_{\di}(1,l)} \\
&+\mathop{\sum \sum}_{m,n \geq 1} \sum_{l \geq 1} \frac{\lambda_{g}(n)}{m \left( nl\right)^{1/2} }  \sum_{\di} \frac{h_{2}(\bm\mu_{\di})}{\mathcal{N}_{\di}} A_{\di}(m,n)  A_{\di}(1,l)\\  
&+\mathop{\sum \sum}_{m,n \geq 1} \sum_{l \geq 1} \frac{\lambda_{g}(n)}{m \left( nl\right)^{1/2} }  \sum_{\di} \frac{h_{3}(\bm\mu_{\di})}{\mathcal{N}_{\di}}\overline{ A_{\di}(m,n)}  A_{\di}(1,l)\\
&+\mathop{\sum \sum}_{m,n \geq 1} \sum_{l \geq 1} \frac{\lambda_{g}(n)}{m \left( nl\right)^{1/2} }  \sum_{\di} \frac{h_{4}(\bm\mu_{\di})}{\mathcal{N}_{\di}}\ovl{ A_{\di}(m,n)}\,  \ovl{A_{\di}(1,l)}\\
&=\mathcal{S}_{1}(T) + \mathcal{S}_{2}(T) + \mathcal{S}_{3}(T) + \mathcal{S}_{4}(T),\,\, \text{(say)}.
 \end{align*}
Here
\begin{align*} h_{1}(\bm\mu_{\di}) = h(\bm\mu_{\di}) W_{\di}(nm^2) \tilde{V}_{\di}(l),\;\;
h_{2}(\bm\mu_{\di}) = h(\bm\mu_{\di})  W_{\di}(nm^2) V_{\di}(l) ,\\
 h_{3}(\bm\mu_{\di}) = h(\bm\mu_{\di})\tilde{W}_{\di}(nm^2)  V_{\di}(l)
 \; \text{and}\;
h_{4}(\bm\mu_{\di}) = h(\bm\mu_{\di}) \tilde{W}_{\di}(nm^2)  \tilde{V}_{\di}(l).
\end{align*}
 Now, Theorem \ref{Theo} follows immediately by the proposition given below.
\begin{proposition}\label{prop} We have
	\begin{align*} 
	\mathcal{S}_{1}(T)&=L(1,g)\frac{1}{192\,\pi^{5}}\iint_{\Re(\bm\mu)=0}h(\bm{\mu})\prod_{j=1}^{3}\frac{\Gamma(\frac{1}{4}+\frac{\mu_{j}}{2})}{\Gamma(\frac{1}{4}-\frac{\mu_{j}}{2})}\mathrm{spec}(\bm\mu)\rm d\bm\mu  + O(T^{\frac{17}{6}+\ve} R^{2}) ,\\	  	  
	\mathcal{S}_{2}(T)&=\zeta(\tfrac{3}{2})\frac{1}{192\pi^{5}}\iint_{\Re(\bm\mu)=0}h(\bm{\mu})\mathrm{spec}(\bm\mu)\rm d\bmu + O(T^{\frac{17}{6}+\ve} R^{2}),\\
	\mathcal{S}_{3}(T)&=L(1,g)\frac{1}{192\,\pi^{5}}\iint_{\Re(\bm\mu)=0}h(\bm{\mu})\prod_{j=1}^{3}\frac{\Gamma(\frac{k}{2}+\mu_{j})}{\Gamma(\frac{k}{2}-\mu_{j})}\mathrm{spec}(\bm\mu)\rm d\bm\mu  + O(T^{\frac{17}{6}+\ve} R^{2}),\\
	\mathcal{S}_{4}(T)&=\zeta(\tfrac{3}{2})\frac{1}{192\pi^{5}}\iint_{\Re(\bmu)=0}h(\bm{\mu})\prod_{j=1}^{3}\frac{\Gamma(\frac{k}{2}+\mu_{j})\Gamma(\frac{1}{4}+\frac{\mu_{j}}{2})}{\Gamma(\frac{k}{2}-\mu_{j})\Gamma(\frac{1}{4}-\frac{\mu_{j}}{2})}\rm{spec}(\bm\mu)\rm d\bmu + O(T^{\frac{17}{6}+\ve} R^{2}).	  	  
	\end{align*}
\end{proposition}

\subsection{Proof of the Proposition \ref{prop}} We only prove the first identity ($S_{1}(T)$) of the Proposition \ref{prop}, as the proof of other identities is the same as the first identity.

 Applying the Kuznetsov's trace formula (Lemma \ref{Kuz}), one has 
 \begin{align}\label{S1}
   S_{1}(T)=\mathop{\sum \sum}_{m,n \geq 1} \sum_{l \geq 1} \frac{\lambda_{g}(n)}{m \left( nl\right)^{1/2} }\left(\Delta^{(1)} + \Sigma_{4}^{(1)} + \Sigma_{5}^{(1)} + \Sigma_{6}^{(1)} -\mathcal{E}_{\rm min}^{(1)} -\mathcal{E}_{ \rm max}^{(1)}\right),
 \end{align}
 where
 \begin{align*}
 \Delta^{(1)}&:=\delta_{m,1}\,\delta_{n,l}\frac{1}{192\pi^{5}}\iint_{\Re(\bmu)=0}h_{1}(\bmu)\rm d_{spec}\bmu,\\
 \Sigma_{4}^{(1)}&:=\sum_{\ep=\pm 1}\sum_{\substack{D_2 \mid D_1 \\  lD_1 =m {D_2}^2}}\frac{\tilde{S}(-\ep n , l , 1 ; D_2 , D_1)}{D_1 D_2} \Phi_{w_4}^{(1)}\left(\frac{\ep ln }{D_1 D_2}\right),\\
 \Sigma_{5}^{(1)} &:=\sum_{\ep=\pm 1}\sum_{\substack{D_1 \mid D_2 \\  D_2 =n {D_1}^2}}\frac{\tilde{S}(\ep m , 1 , l ; D_1 , D_2)}{D_1 D_2} \Phi_{w_5}^{(1)}\left(\frac{\ep lm}{D_1 D_2}\right),\\
 \Sigma_{6}^{(1)} &:=\sum_{\ep_1 ,\ep_2 =\pm 1}\sum_{D_1 ,D_2}\frac{S(\ep_2 n ,\ep_1 m , 1 , l ; D_1 , D_2)}{D_1 D_2} \Phi_{w_6}^{(1)}\left(-\frac{\ep_2 n  D_2}{{D_1}^{2}},-\frac{\ep_1lm  D_1}{{D_2}^{2}} \right), 
 \end{align*}
with $\Phi_{w_4}^{(1)}(y)$, $\Phi_{w_5}^{(1)}(y)$ and $\Phi_{w_6}^{(1)}(\bm y)$ defined as in \eqref{Phi} by using the new test function $h_{1}(\bm\mu_{\di}) = h(\bm\mu_{\di}) W_{\di}(nm^2) \tilde{V}_{\di}(l)$ respectively; and 

\begin{align*}
\mathcal{E}_{\rm min}^{(1)}&:=\frac{1}{24 (2 \pi i)^2} {\iint}_{\Re(\bmu)=0} \frac{h_{1}(\bmu)}{\mathcal{N}_{\bmu}^{\rm min}} \ovl{B_{\bmu}^{\min}(1,l)} B_{\bmu}^{\min}(m,n) \rm d\mu_1 d\mu_2,\\
\mathcal{E}_{\rm max}^{(1)}&:= \sum_{f}\frac{1}{2 \pi i} \int_{\Re(u)=0}\frac{h_{1}(u+it_f , u-it_f , -2u )}{\mathcal{N}_{u,f}^{\rm max}}\ovl{B_{u,f}^{\rm max}(1 ,l)} B_{u,f}^{\rm max}(m ,n) {\rm d}u,
\end{align*}
\subsubsection{\bf The diagonal term} Let us denote $\mathcal{D}^{(1)}$ to be contribution of $\Delta^{(1)}$ to the sum in the equation \eqref{S1}. Thus we have 
\begin{align*}
\mathcal{D}^{(1)}= \frac{1}{192\pi^{5}}\iint_{\Re(\bmu)=0} h(\bmu) D(\bmu) \rm d_{spec}\bmu,
\end{align*}
 where
 \begin{align} \label{eq6}
 D(\bmu):=\frac{1}{(2\pi i)^{2}}\int_{(3)}\int_{(3)}L(1+u_{1} +u_{2}, g )
 \frac{\gamma(\tfrac{1}{2}+u_1 ,\tilde{\di})}{\gamma(\tfrac{1}{2} ,\di)}\frac{\gamma(\tfrac{1}{2}+u_2 ,g\ot\di)}{\gamma(\tfrac{1}{2}, g\ot\di)} \nonumber \\
 \times G(u_1)G(u_2)\frac{{\rm d}u_1}{u_1}\frac{{\rm d}u_2}{u_2} 
 \end{align}
with $G(u)=e^{u^2}$, $\gamma(u,\tilde{\di})=\displaystyle\prod_{j=1}^{3}\Gamma_{\mathbb{R}}(u+\mu_j)$,  $\gamma(u,\di)=\displaystyle\prod_{j=1}^{3}\Gamma_{\mathbb{R}}(u-\mu_j)$ and
$$  \gamma(u,g\ot \di)=\prod_{j=1}^{3}\Gamma_{\mathbb{R}}(u+\tfrac{k-1}{2}-\mu_j) \Gamma_{\mathbb{R}}(u+\tfrac{k+1}{2}-\mu_j).$$ 
Now we evaluate the double integral in the first term on the right side of the equation \eqref{eq6}. We first shift the contour to ${\Re}(u_1)=\ep $, ${\Re}(u_2)=\ep$ without encountering a pole, for any $\ep>0$. Then we move the contour ${\Re}(u_2)=\ep$ to ${\Re}(u_2)=-\frac{1}{2}$, in doing so we encounter a simple pole at $u_2 =0$. Since, $G(u)\ll e^{-t^2}$ and $L(1/2 +it, g)\ll t^{1/3}$, the integral on ${\Re}(u_1)=\ep$, ${\Re}(u_2)=-\frac{1}{2}$ is bounded by $O\left(\textstyle\prod_{j=1}^{3}|\mu_{j}|^{-1/2}\right)$. The contribution from the residue at $u_2 =0$ is 
\begin{align}\label{eq7}
\frac{1}{2\pi i}\int_{(\epsilon)}L(1+u_{1} ,g)
\frac{\gamma(\tfrac{1}{2}+u_1 ,\tilde{\di})}{\gamma(\tfrac{1}{2} ,\di)}G(u_1)\frac{{\rm d}u_1}{u_1}.
\end{align}
Now shift the contour in \eqref{eq7} to ${\Re}(u_1)=-\frac{1}{2}$ encountering a simple pole at $u_1 =0$. The integral on the line ${\Re}(u_1)=-\frac{1}{2} $ is bounded by $O\left(\textstyle\prod_{j=1}^{3}|\mu_{j}|^{-1/4}\right)$. The contribution from the residue at $u_1 =0$ is $L(1,g)\displaystyle\frac{\gamma(\frac{1}{2}, \tilde{\di})}{\gamma(\frac{1}{2}, \di)}$.
Note that, $\displaystyle\iint_{\Re(\bmu)=0} h(\bm{\mu})\mathrm{spec}(\bm\mu)d\bm\mu\asymp T^{3}R^{2}$. Therefore, 
\[
\mathcal{D}^{(1)}=L(1,g)\frac{1}{192\pi^{5}}\iint_{\Re(\bm\mu)=0} h(\bm{\mu})\prod_{j=1}^{3}\frac{\Gamma(\frac{1}{4}+\frac{\mu_{j}}{2})}{\Gamma(\frac{1}{4}-\frac{\mu_{j}}{2})}\mathrm{spec}(\bm\mu)\rm d\bmu + O(T^{\frac{9}{4}+\ve} R^{2}).	 
\]  
 
\subsubsection{\bf Contribution of $\Sigma_{4}^{(1)}$} Let $E^{(1)}_{4}$ be the contribution of $\Sigma_{4}$ to the sum in the equation \eqref{S1}. So
\[
E^{(1)}_{4}=\mathop{\sum \sum}_{m,n \geq 1} \sum_{l \geq 1}\frac{\lambda_{g}(n)}{m \left( nl\right)^{1/2} }\sum_{\ep=\pm 1}\sum_{\substack{D_2 \mid D_1 \\  lD_1 =m {D_2}^2}}\frac{\tilde{S}(-\ep n , l , 1 ; D_2 , D_1)}{D_1 D_2} \Phi_{w_4}^{(1)}\left(\frac{\ep ln }{D_1 D_2}\right).
\] 
Let $U_{i}(x)$ $(i=1,2)$ be smooth functions which are compactly supported in $[1, 2]$, satisfy $x^{j}U_{i}^{(j)}(x)\ll 1$. By partition of unity, to get a bound for $E^{(1)}_{4}$ it is enough to estimate the following sum
\begin{align}\label{eq8}
\mathop{\sum \sum}_{m,n \geq 1} \sum_{l \geq 1}\frac{\lambda_{g}(n)U_{1}(\tfrac{nm^2}{N})U_{2}(\frac{l}{L})  }{ \left(nm^2\right)^{1/2}(l)^{1/2} } \sum_{\substack{\de , D \\  \de l=m D}} \frac{\tilde{S}(\mp n , l , 1 ; D , \de D)}{\de D^2} \Phi_{w_4}^{(1)}\left(\frac{\pm ln }{\de D^2}\right),
\end{align} where $1\le N\le T^{3+\ve}$ and $1\le L\le T^{3/2\,+\ve}$. By Lemma \ref{le4}, $\Phi_{w_4}^{(1)}\left(\frac{\pm ln }{\de D^2}\right)$ is negligibly small unless
\begin{align}\label{eq9}
\frac{ ln }{\de D^2}\gg T^{3-\ve}.
\end{align}
Note that $\de l=mD$. By \eqref{eq9}, we deduce that 
\begin{align*}
1\le l\de^3\le\frac{nm^2}{T^{3-\ve}}\le T^{\ve},
\end{align*} 
which implies, $\de, L, mD\le T^{\ve}$ and $T^{3-\ve}\le n\le N\le T^{3+\ve}$.

\noindent
We recall the Kloosterman sum $$\tilde{S}(\mp n , l , 1 ; D , \de D)=\mathop{\sum \sum}_{\substack{C_1({\rm mod}\, D ),\, C_2({\rm mod}\, \de D ) \\ (C_1 ,D)=(C_2 ,\de)=1 }} e\left(l \frac{\bar{C_{1}} C_2}{ D} +  \frac{\bar{C_{2}}}{\de} \mp  n \frac{C_1}{D} \right). $$
In \eqref{eq8}, the only non-trivial sum is the sum over $n$, which is given by 
\[
\sum_{n \geq 1}\lambda_{g}(n)\,e\left(\mp n \frac{C_1}{D}\right) \theta(n),
\] 
where 
\[ 
\theta(y)=\frac{1}{\sqrt{y}}U_{1}\left(\frac{m^2 y}{N}\right)\Phi_{w_4}^{(1)}\left(\frac{\pm ly }{\de D^2}\right).
\] 
Now the $GL(2)$-Voronoi summation formula transforms the above sum into
\[
\frac{1}{D}\sum_{n \geq 1}\lambda_{g}(n)\, e\left(\pm n \frac{\bar{C_1}}{D}\right)\Theta(n),
\]
where 
\[
\Theta(y)=2\pi i^{k}\int_{0}^{\infty}\theta(x)J_{k-1}\left(\frac{4\pi\sqrt{xy}}{D}\right){\rm d}x.
\] 
We now analyse the integral transform $\Theta(n)$. Using the properties of Bessel function we arrive at the following expression
\begin{align}\label{eq10}
\frac{N^{1/4} D^{1/2}}{m^{1/2}n^{1/4}}\int_{0}^{\infty}\frac{1}{x^{3/4}} U_{1}(x) \Phi_{w_4}^{(1)}\left(\frac{\pm lNx }{\de (m D)^2}\right)\, e\left(\pm\frac{2\pi \sqrt{nNx}}{mD} \right){\rm d}x.
\end{align} 
By Lemma \ref{le4}, we have 
\[
|x|^{j}\frac{\rm d^{j}\,}{\;{\rm d}x^{j}}\Phi_{w_4}^{(1)}(x)\ll_{\ve,j} T^{3+\ve}R^{2}\left(T+|x|^{1/3}\right)^{j}.
\]
First, we change the variable $x\to x^2$ in the integral of equation \eqref{eq10}. Then by repeated integration by parts we see that
\[
\Theta(n)\ll_{\ve, j}\frac{N^{1/4} D^{1/2}T^{3+\ve}R^{2}}{m^{1/2}n^{1/4}}\left(T+\left(\frac{lN}{\de (mD)^2}\right)^{1/3}\right)^{j}\left(\frac{mD}{\sqrt{nN}}\right)^{j}.
\] 
Since $1\le\de, L, m,D\le T^{\ve}$ and $T^{3-\ve}\le N\le T^{3+\ve}$.

\noindent
Therefore, $$E_{4}^{(1)}\ll_{A} T^{-A},$$
for any $A>0$.
 
\subsubsection{\bf Contribution of $\Sigma_{5}^{(1)}$} Let $E^{(1)}_{5}$ be the contribution of $\Sigma_{5}$ to the sum in the equation \eqref{S1}. As in the previous case it is enough to estimate the following sum 
\begin{align}\label{eq11}
\mathop{\sum \sum}_{m,n \geq 1} \sum_{l \geq 1}\frac{\lambda_{g}(n)U_{1}(\tfrac{nm^2}{N})U_{2}(\frac{l}{L})}{\left(nm^2\right)^{1/2}(l)^{1/2}} \sum_{\substack{\de , D \\  \de =n D}}\frac{\tilde{S}(\pm m , 1 , l ; D ,\de D)}{\de D^2} \Phi_{w_5}^{(1)}\left(\frac{\pm lm}{\de D^2}\right)
\end{align}
to get a bound for $E^{(1)}_{5}$. Here $1\le nm^2 \le N\le T^{3+\ve}$ and $1\le l\le L\le T^{3/2\,+\ve  }$. Since $\Phi_{w_5}^{(1)}\left(\frac{\pm lm}{\de D^2}\right)$ is negligibly small unless
\begin{align*}
T^{3-\ve}\ll \frac{lm}{\de D^2}.
\end{align*}
Note that $\de=n D$, together with above inequality we get 
\[
nD^{3}\le \frac{lm}{T^{3-\ve}}\le T^{\ve}.
\] 
Therefore $n,\de, D\le T^{\ve}$ and $T^{3/2\,-\ve}\le m,l\le T^{3/2\,+\ve}$. 
Recall the Kloosterman sum $$\tilde{S}(\pm m , 1 , l ; D , \de D)=\mathop{\sum \sum}_{\substack{C_1({\rm mod}\, D ),\, C_2({\rm mod}\, \de D ) \\ (C_1 ,D)=(C_2 ,\de)=1 }} e\left( \frac{\bar{C_{1}} C_2}{ D} +  l \frac{\bar{C_{2}}}{\de} \mp  m \frac{C_1}{D} \right). $$ Now, the $l$-sum in the equation \eqref{eq11} is given by
\[
\sum_{l\ge 1}e\left(\frac{l\bar{C_{2}}}{\de}\right)\frac{1}{\sqrt{l}}U_{2}\left(\frac{l}{L}\right)\Phi_{w_5}^{(1)}\left(\frac{\pm lm}{\de D^2}\right).
\]  
An application of the Poisson summation formula to the above sum yields
\begin{align}\label{eq12}
\sqrt{L}\sum_{\substack{l\in \mathbb{Z}\\ l\equiv-C_{2} (\rm mod\; \de)}}\int_{\mathbb{R}}\frac{1}{\sqrt{y}}U_{2}(y)\Phi_{w_5}^{(1)}\left(\frac{\pm Lmy}{\de D^2}\right) e\left(\frac{-lLy}{\de}\right){\rm d}y.
\end{align}
Note that, for $l\neq 0$ (non-zero frequency), by repeated integration by parts we have the inner integral in \eqref{eq12} is 
\[
\ll_{\ve, j} T^{3+\ve}R^{2}\left(T+\left(\frac{mL}{\de D^{2}}\right)^{1/3}\right)^{j} \left(\frac{\de}{|l|L}\right)^{j}.
\]
Thus we need to have $\de=1$, otherwise there will be no zero frequency. Therefore, we have $D=1$ and $n=1$. In this case the contribution of zero frequency is given by 
\begin{align*}
\sqrt{L}\sum_{m \geq 1}\frac{1}{m}U_{1}\left(\frac{m^2}{N}\right)\int_{\mathbb{R}}\frac{1}{\sqrt{y}}U_{2}(y) \Phi_{w_5}^{(1)}\left(\pm Lmy\right){\rm d}y.
\end{align*}
By the definition of $\Phi_{w_5}^{(1)}(y)$ and $K_{w_4}(y;\bmu)$, the above $y$-integral becomes
\begin{align}\label{eq13}
&\int_{\mathbb{R}}\frac{1}{\sqrt{y}}U_{2}(y)\iint_{\Re(\bmu)=0}h_{1}(\bmu)K_{w_4}\left(\mp Lm y;-\bmu\right) {\rm d_{spec}}\bmu\,{\rm d}y\nonumber\\
=&\iint_{\Re(\bmu)=0}h_{1}(\bmu)\int_{-i \infty}^{i\infty}\left(\int_{\mathbb{R}}y^{-s-\frac{1}{2}} U_{2}(y){\rm d}y\right)|mL|^{-s}\tilde{G}^{\pm}(s,-\bmu)\frac{{\rm d}s}{2\pi i}\, {\rm d_{spec}}\bmu\nonumber\\
=&\iint_{\Re(\bmu)=0}h_{1}(\bmu)\int_{-i \infty}^{i\infty} |mL|^{-s} \widehat{U}_{1}\left(\frac{1}{2}-s\right)\tilde{G}^{\pm}(s,-\bmu)\frac{{\rm d}s}{2\pi i}\, {\rm d_{spec}}\bmu,
\end{align}
where $\widehat{U}_{1}(s)$ is the Mellin transform of $U_{1}(y)$, which is entire and rapidly decaying. Then we can restrict the $s$-integral to $|{\rm Im}(s)|\le T^{\ve}$. Recall the definition of $\tilde{G}^{\pm}(s,\bmu)$
\[  
\tilde{G}^{\pm}(s,\bmu):=\frac{\pi^{-3s}}{12288\pi^{7/2}}\left(\prod_{j=1}^{3}\frac{\Gamma(\tfrac{1}{2}(s-\mu_j))}{\Gamma(\frac{1}{2}(1-s+\mu_j))}\pm i \prod_{j=1}^{3}\frac{\Gamma(\frac{1}{2}(1+s-\mu_j))}{\Gamma(\frac{1}{2}(2-s+\mu_j))} \right).
\]
Since $\bmu_0 \asymp\mu_{0,j} \asymp T$, then the $\bmu$-integral in \eqref{eq13} is bounded by $T^{3/2\, +\ve}R^2$.

Thus, $$E^{(1)}_{5}\ll T^{\frac{9}{4}+\ve}R^2.$$ 
 
\subsubsection{\bf Contribution of $\Sigma_{6}^{(1)}$}
In this case we have to bound the following sum
\begin{align}\label{eq14}
\mathop{\sum \sum}_{m,n \geq 1} \sum_{l \geq 1} \frac{\lambda_{g}(n)U_{1}(\tfrac{nm^2}{N})U_{2}(\frac{l}{L})}{\left(nm^2\right)^{1/2}(l)^{1/2}} \sum_{D_1 ,D_2}\frac{S(\pm n ,\pm m , 1 , l ; D_1 , D_2)}{D_1 D_2} \Phi_{w_6}^{(1)}\left(\frac{\mp n  D_2}{{D_1}^{2}},\frac{\mp lm  D_1}{{D_2}^{2}}\right).
\end{align} 
Using the property $\Phi_{w_6}$ (Lemma \ref{le5}), we have
\[
T^{1-\ve}\le \frac{(l m n^2)^{1/6}}{D_{1}^{1/2}}\quad\text{and}\quad T^{1-\ve}\le \frac{(n l^2 m^2)^{1/6}}{D_{2}^{1/2}}.
\] 
From these conditions we infer that $1 \leq D_{2} \leq T^{\ve}, 1 \leq D_{1} \leq T^{1/2+\ve}/ m$ and $1 \leq m \leq T^{1/2+\ve}$, also $ T^{3-\ve} \leq N \leq T^{3+\ve}$ and $T^{3/2 -\ve} \leq L \leq T^{3/2+\ve}$.
 
\noindent
 The kloosterman sun $S(\pm n ,\pm m , 1 , l ; D_1 , D_2 )$ is given by
\begin{align*}
\mathop{\sum \sum \sum \sum}_{\substack{B_1 , C_1({\rm mod}\, D_1 );\, B_2 , C_2({\rm mod}\, D_2 )\\ D_1 C_2 + B_1 B_2 + D_2 C_1 \equiv 0\, ({\rm mod}\, D_1 D_2 ) \\ (B_j  ,C_j ,D_j)=1 }} e\left( \frac{\pm n B_1 +  (Y_1 D_2 - Z_1 B_2)}{D_1} + \frac{\pm m B_2 + l(Y_2 D_1 - Z_2 B_1)}{D_2} \right),
\end{align*}
 where $  B_j Y_j + C_j  Z_j\equiv 1\, ({\rm mod}\, D_j )$ for $j=1,2$.
Note that $(B_{1}, D_{1})\mid D_2$, so $(B_{1}, D_{1})\le T^{\ve}$ as $D_2 \le T^{\ve}$. Let $B_{1}=B^{'}_{1} (B_{1}, D_{1})$ and $D_{1}=D^{'}_{1} (B_{1}, D_{1})$  with $(B^{'}_{1},D^{'}_{1})=1$. 

\noindent
Now we apply $GL(2)$-Voronoi summation and Poisson summation formulae on $n$ and $l$-sums in the equation \eqref{eq14} respectively. Then the $n$ and $l$-sums in \eqref{eq14} transform into
\begin{align}\label{eq15}
&2\pi i^k \frac{m\sqrt{L}}{\sqrt{N}D^{'}_{1}}\sum_{n \geq 1}\sum_{\substack{l\in\mathbb{Z}\\l\equiv-(Y_2 D_1 - Z_2 B_1)\,{\rm mod}\,D_{2}}}\lambda_{g}(n)e\left(\frac{\mp \bar{B}^{'}_{1}n}{D^{'}_{1}}\right)\nonumber\\ &\times\int_{\mathbb{R}}\int_{\mathbb{R}}\Phi_{w_6}^{(1)}\left(\frac{\mp N D_2 x}{{D_1}^{2} m^2},\frac{\mp m L  D_1 y}{{D_2}^{2}}\right) J_{k-1}\left(\frac{4\pi \sqrt{Nnx}}{mD_{1}^{'}}\right)e\left(\frac{-lLy}{D_{2}}\right)U_{1}(x)U_{2}(y){\rm d}x\, {\rm d}y.
\end{align}
By repeated integration by parts, the $y$-integral in \eqref{eq15} is 
\[
\ll_{\ve, j} T^{3+\ve}R^{2}\left(T+\left(\frac{mD_{1}L}{D_{2}^{2}}\right)^{1/2}+\left(\frac{mD_{1}L}{D_{2}^{2}}\right)^{1/3}\left(\frac{ND_{2}}{D_{1}^{2}m^{2}}\right)^{1/6}\right)^{j}\left(\frac{D_2}{|l|L}\right)^{j},
\]  
for $l\neq 0$. Therefore the non-zero frequency ($l\neq 0$) in \eqref{eq15} contributes $O(T^{-A})$. In the zero frequency $(l= 0)$ case, we must have $D_{2}\mid(Y_2 D_1 - Z_2 B_1)$.  
In this case equation \eqref{eq15} becomes
\begin{align*}
2\pi i^k \frac{m\sqrt{L}}{\sqrt{N}D^{'}_{1}}\sum_{n \geq 1}\lambda_{g}(n)e\left(\frac{\mp \bar{B}^{'}_{1}n}{D^{'}_{1}}\right) \int_{\mathbb{R}}\int_{\mathbb{R}}\Phi_{w_6}^{(1)}\left(\frac{\mp N D_2 x}{{D_1}^{2} m^2},\frac{\mp m L  D_1 y}{{D_2}^{2}}\right)\\
 \times\, J_{k-1}\left(\frac{4\pi \sqrt{Nnx}}{mD_{1}^{'}}\right)U_{1}(x)U_{2}(y){\rm d}x\, {\rm d}y.
\end{align*}
Inserting the properties of Bessel function in the above equation, we have 
\begin{align}\label{eq16}s
\frac{m^{3/2}\sqrt{L}}{N^{3/4}{D^{'}_{1}}^{1/2}}\sum_{n \geq 1}\frac{\lambda_{g}(n)}{n^{1/4}}e\left(\frac{\mp \bar{B}^{'}_{1}n}{D^{'}_{1}}\right)  \int_{\mathbb{R}}\int_{\mathbb{R}}\Phi_{w_6}^{(1)}\left(\frac{\mp N D_2 x}{{D_1}^{2} m^2},\frac{\mp m L  D_1 y}{{D_2}^{2}}\right)\nonumber\\
\times\, e\left(\frac{2\pi \sqrt{Nnx}}{mD_{1}^{'}}\right)U_{1}(x)U_{2}(y){\rm d}x\, {\rm d}y.
\end{align}
Integrating by parts we see that the $x$-integral in \eqref{eq16} is negligibly small unless
$n\le T^{\ve}$ (Dual length). Therefore the sum in \eqref{eq14} is bounded above by
\[
\frac{T^{\ve}\sqrt{L}}{N^{3/4}}\sum_{n \le T^{\ve}}\frac{\lambda_{g}(n)}{n^{1/4}}\mathop{\sum\sum}_{\substack{D_{1}\le T^{1/2 +\ve}\\ D_{2}\le T^{\ve}}} \frac{1}{D_{1}^{3/2}D_{2}}\sum_{m\le \frac{T^{1/2 +\ve}}{D_{1}}}m^{1/2} \mathcal{C}(D_{1},D_{2},n,m)\mathcal{I}(D_{1},D_{2},n,m),
\]
where \begin{align*}
\mathcal{C}(D_{1},D_{2},n,m)&=\mathop{\sum\sum}_{\substack{B_1 , C_1({\rm mod}\, D_1 );\, B_2 , C_2({\rm mod}\, D_2 )\\ D_1 C_2 + B_1 B_2 + D_2 C_1 \equiv 0\, ({\rm mod}\, D_1 D_2 ) \\ (B_j  ,C_j ,D_j)=1 \\ D_{2}\mid (Y_2 D_1 - Z_2 B_1)}} e\left(\frac{\mp \bar{B}^{'}_{1}n }{D^{'}_{1}}+\frac{Y_1 D_2 - Z_1 B_2}{D_1}\pm \frac{mB_{2}}{D_2}\right),    \\
\mathcal{I}(D_{1},D_{2},n,m)&=\int_{\mathbb{R}}\int_{\mathbb{R}}\Phi_{w_6}^{(1)}\left(\frac{\mp N D_2 x}{{D_1}^{2} m^2},\frac{\mp m L  D_1 y}{{D_2}^{2}}\right) e\left(\frac{2\pi \sqrt{Nnx}}{mD_{1}^{'}}\right)U_{1}(x)U_{2}(y){\rm d}x\, {\rm d}y.
\end{align*} 
Moreover, trivially $\mathcal{C}(D_{1},D_{2},n,m)$ is bounded by $O(T^{\ve} D_1)$ and the integral $$\mathcal{I}(D_{1},D_{2},n,m)\ll T^{3+\ve}R^2.$$
 Hence, the contribution of $\Sigma_{6}^{(1)}$ is $O(T^{\frac{9}{4}+\ve}R^2)$.  

\subsubsection{\bf Contribution of $\mathcal{E}_{\min}^{(1)}$} Let us denote
\begin{align*}
E^{(1)}_{\min}:=\frac{1}{24 (2 \pi i)^2}{\iint}_{\Re(\bmu)=0}\; \frac{h(\bmu)}{\mathcal{N}_{\bmu}^{\rm min}}\mathop{\sum \sum}_{m,n \geq 1} \frac{\lambda_{g}(n)}{m n^{1/2}}B_{\bmu}^{\min}(m,n)W_{\di}(nm^2)\\
\times \sum_{l \geq 1}\frac{1}{l^{1/2}}\ovl{B_{\bmu}^{\min}(1,l)}\tilde{V}_{\di}(l)\, \rm d\mu_1 d\mu_2.
\end{align*}
Inserting the definition of $V_{\di}(y)$, $\tilde{W}(y)$ and using \eqref{gEmin}, \eqref{Emin} we have
\[
E^{(1)}_{\min}=\frac{1}{24 (2 \pi i)^2}{\iint}_{\Re(\bmu)=0}\; \frac{h(\bmu)}{\mathcal{N}_{\bmu}^{\rm min}}\mathcal{I}^{(1)}_{\min}(\bmu)\rm d\mu_1 d\mu_2,
\] where
\begin{align*}
\mathcal{I}&^{(1)}_{\min}(\bmu)=\frac{1}{(2\pi i)^2} \int_{(3)}\int_{(3)}  
G(s_2)\prod_{j=1}^{3}\frac{\Gamma_{\mathbb{R}} \left(s_{2}+\frac{1}{2} +\mu_{j} \right)}{\Gamma_{\mathbb{R}} \left(\frac{1}{2}-\mu_{j} \right)}\zeta(s_2 +\tfrac{1}{2}-\mu_{j})
\\
&\times G(s_1)\prod_{j=1}^{3}\frac{\Gamma_{\mathbb{R}} \left(s_{1}+\frac{1}{2}+ \frac{k-1}{2}-\mu_{j} \right)  \Gamma_{\mathbb{R}} \left(s_{1}+\frac{1}{2}+ \frac{k+1}{2}-\mu_{j} \right)}{\Gamma_{\mathbb{R}} \left(\frac{1}{2}+ \frac{k-1}{2}-\mu_{j} \right)  \Gamma_{\mathbb{R}} \left(\frac{1}{2}+ \frac{k+1}{2}-\mu_{j} \right)}L(s_{1}+\tfrac{1}{2}+\mu_j , g)\frac{{\rm d}s_1}{s_1}\frac{{\rm d}s_2}{s_2}.
\end{align*}
Now we move the line of integration from $\Re(s_1)=\Re(s_2)=3$ to $\Re(s_1)=\Re(s_2)=\ve$. Since $L(1/2 +it,g)\ll t^{1/3}$, $\zeta(1/2+ it)\ll t^{1/6}$ and $G(s)\ll e^{-t^2}$, therefore 
$$\mathcal{I}^{(1)}_{\min}(\bmu)\ll\prod_{j=1}^{3} (1+|\mu_j |)^{1/2}.$$
Finally, using the bound $$\mathcal{N}_{\bmu}^{\rm min}=\mathcal{N}_{\nu_1 , \nu_2}^{\rm min}:= \frac{1}{16}\prod_{j=1}^{3} |\zeta(1+3\nu_j)|^{2}\gg \prod_{j=1}^{3}\left(\frac{1}{\log(1+3\Im(\nu_j))}\right)^{2},$$ we obtain $${E}_{\min}^{(1)}\ll T^{\frac{3}{2}+\ve}R^2 .$$

\subsubsection{\bf Contribution of $\mathcal{E}_{\max}^{(1)}$}\label{S3.16}

Let us define \begin{align*} E_{\max}^{(1)}:=\sum_{f}\frac{1}{2 \pi i} &\int_{\Re(u)=0}\frac{h(u+it_f , u-it_f , -2u )}{\mathcal{N}_{u,f}^{\rm max}}\\
&\times \mathop{\sum \sum}_{m,n \geq 1} \frac{\lambda_{g}(n)}{m n^{1/2}} B_{u,f}^{\rm max}(m ,n)W_{\di}(nm^2)\sum_{l \geq 1}\frac{1}{l^{1/2}}\ovl{B_{u,f}^{\rm max}(1 ,l)}\tilde{V}_{\di}(l) {\rm d}u.
\end{align*}
By similar argument as above one has
\[
E^{(1)}_{\max}=\sum_{f}\frac{1}{2 \pi i} \int_{\Re(u)=0}\frac{h(u+it_f , u-it_f , -2u )}{\mathcal{N}_{u,f}^{\rm max}}\mathcal{I}(u+it_f , u-it_f , -2u) {\rm d}u,
\] where $\mathcal{I}^{(1)}_{\max}(u+it_f , u-it_f , -2u ) $ is given by
\begin{align*}
&\frac{1}{(2\pi i)^2} \int_{(3)}\int_{(3)} G(s_1)  
G(s_2)L(s_{1}+\tfrac{1}{2}-2u, g) L(s_{1}+\tfrac{1}{2}+u, g\ot f)
L(s_{2}+\tfrac{1}{2}-u, g)\times
\\
&\zeta(s_2 +\tfrac{1}{2}+2u)\prod_{j=1}^{3}\frac{\Gamma_{\mathbb{R}}\left(s_{1}+\frac{1}{2}+ \frac{k-1}{2}-\alpha_{j} \right)  \Gamma_{\mathbb{R}} \left(s_{1}+\frac{1}{2}+ \frac{k+1}{2}-\alpha_{j} \right) \Gamma_{\mathbb{R}} \left(s_{2}+\frac{1}{2} +\alpha_{j} \right)}{\Gamma_{\mathbb{R}} \left(\frac{1}{2}+ \frac{k-1}{2}-\alpha_{j} \right)  \Gamma_{\mathbb{R}} \left(\frac{1}{2}+ \frac{k+1}{2}-\alpha_{j} \right) \Gamma_{\mathbb{R}} \left(\frac{1}{2}-\alpha_{j} \right)}\frac{{\rm d}s_1}{s_1}\frac{{\rm d}s_2}{s_2}.
\end{align*}
Here $\alpha_{1}=u+it_f$, $\alpha_{2}= u-it_f$ and $\alpha_{3}= -2u$. By the definition of $h$ we have $h(u+it_f , u-it_f , -2u )$ is negligibly small unless 
\[
 | u+it_f -\mu_{0,1}|\le R ,\;\; |u-it_f -\mu_{0,2}|\le R,\;\; |-2u-\mu_{0,3}|\le R .
\]
Note that $\mu_{0,j}\asymp T$ and $\mathcal{N}_{u,f}^{\max}:= 8L(1,{\rm Ad}^{2}f)|L(1+3u,f)|^{2} \gg (1+\log|u|)^{-1}$.

We first shift the line of integration to $\Re(s_1 )=\Re(s_2 )=\ve$ and use the bounds $$\zeta(1/2+it) \ll |t|^{1/6},\;\, L(1/2+it,f) \ll |t|^{1/3},\;\, L(1/2+it,g) \ll |t|^{1/3}$$ and $L\left(1/2+it, g \otimes f\right) \ll |t|$ to get 
\[
E^{(1)}_{\max}\ll T^{11/6 +\ve}R\sum_{T-R \leq t_{f} \leq T+R}1.
\]
Using the Weyl law for the number of eigenvalues associated to $GL(2)$ cusp forms we have
$$E^{(1)}_{\max}\ll T^{17/6 +\ve}R^2 .$$

\vspace{1cm}
\noindent
\textbf{Acknowledgements:} We thank Ritabrata Munshi, Satadal Ganguly and D. Surya Ramana for their encouragement and support. We also thank the Stat-Math unit, Indian Statistical Institute Kolkata for the excellent academic atmosphere. The second-named author is supported by the NBHM grant(No: 0204/16(14)/2021/R\&D-II/01).

\end{document}